
\magnification=\magstep1
\baselineskip=14pt
\def \Box {\hbox{}\nobreak
 \vrule width 1.6mm height 1.6mm depth 0mm  \par \goodbreak \smallskip}

\centerline{\bf THE GEOMETRY OF A-GRADED ALGEBRAS}

\vskip .3cm

\centerline{{\bf Bernd Sturmfels}}

\vskip .1cm

\centerline{Courant Institute, New York University, New York, NY 10012, USA}
\centerline{and \ \ University of California, Berkeley, CA 94720, USA}

\vskip .4cm

\midinsert  \narrower {\sl
\baselineskip=11pt
\centerline{Abstract}
\noindent
We study algebras $k[x_1,...,x_n]/I$ which admit a grading by a subsemigroup
of ${\bf N^d}$  such that every graded component is a one-dimensional
$k$-vector space. V.I.~Arnold and coworkers proved that for $d = 1 $
and $n \leq 3$ there are only finitely many isomorphism types of such
${\cal A}$-graded algebras, and in these cases $I$  is an initial ideal
(in the sense of Gr\"obner bases) of a toric ideal.
In this paper it is shown that Arnold's finiteness theorem does not
extend to  $n = 4$. Geometric conditions are given for $I$ to be an initial
ideal of a toric ideal. The varieties defined by ${\cal A}$-graded algebras
are characterized in terms of polyhedral subdivisions, and the distinct
${\cal A}$-graded algebras are parametrized by a certain binomial scheme.
}

\endinsert

\vskip .4cm

\beginsection 1. Introduction

What are the graded algebras that have the simplest possible Hilbert function~?
This question was raised and partially answered by V.I.~Arnold [1] and
his coworkers  E.~Korkina, G.~Post and M.~Roelofs [9],[10],[11].
They considered finitely generated ${\bf N}$-graded $k$-algebras such that each
non-trivial graded component has dimension $1$ over the ground field $k$.
We propose the following multigraded version of their definition:
Let ${\cal A} = \{{\bf a}_1,{\bf a}_2,\ldots,{\bf a}_n \}$ be a
subset of ${\bf N}^d \!\setminus \! \{ 0 \}$,
and let ${\bf N} {\cal A}$ denoted the sub-semigroup of ${\bf N}^d$
spanned by ${\cal A}$. An {\it ${\cal A}$-graded algebra} is a
${\bf N}^d$-graded $k$-algebra $\, R = \bigoplus_{\bf b} R_{\bf b} \,$ with
homogeneous generators  $X_1,X_2,\ldots,X_n$  in degrees ${\bf a}_1,
{\bf a}_2,\ldots, {\bf a}_n$ such that
$$ dim_k (R_{\bf b}) \quad = \quad
\cases{ 1 & if $ {\bf b} \in {\bf N} {\cal A} $ \cr
        0 & otherwise }
\qquad\qquad \hbox{for all } \, {\bf b} \in {\bf N}^d.  \eqno (1.1) $$
Every ${\cal A}$-graded algebra has a natural presentation
as a quotient of a polynomial ring:
$$ 0 \quad \rightarrow \quad I \quad \rightarrow \quad
k[x_1,x_2,\ldots,x_n]  \quad \rightarrow \quad R
 \quad \rightarrow \quad 0 .$$
The presentation ideal $I = ker(x_i \mapsto X_i)$ is
called {\it ${\cal A}$-graded} as well.
It is easy to see  (cf.~[5, Proposition 1.11])
that $I$ is generated by polynomials with at most two terms,
that is, $I$ is a {\it binomial ideal}.
In the following we use the abbreviation
$k[{\bf x}] := k[x_1,x_2,\ldots,x_n]$.

Arnold, Korkina, Post and Roelofs studied the case $d=1$, and
they proved that for $n \leq 3$ there is only a finite number
of non-isomorphic ${\cal A}$-graded algebras.
To state their main result (Theorem 1.1),
we need some definitions.
Two ${\cal A}$-graded algebras $R$ and $R'$ are {\it isomorphic}
if there exists an algebra isomorphism of degree $0$.
This holds if and only if, for the corresponding ideals $I$ and $ I'$,
there exists $\lambda = (\lambda_1,\ldots,\lambda_n) $ in $ (k^*)^n$ such that
$$   I' \quad = \quad \lambda \cdot I \,\,\, := \,\,\
\{ \, f(\lambda_1 x_1, \ldots,
\lambda_n x_n) \,\,: \,\,  f \in I \,\}. \eqno (1.2) $$
Gr\"obner basis theory suggests taking limits with respect
to one-parameter subgroups $\omega$ of the torus $(k^*)^n$.
Each such limit is a  new ${\cal A}$-graded ideal:
the {\it initial ideal} $ in_\omega(I) $, which is spanned by the
highest forms $in_\omega(f)$ of all polynomials $f$ in $I$.
Here $\omega$ is identified with a vector in ${\bf Z}^n$,
and, if $in_\omega(I)$ is a monomial ideal, then
$\omega $ is called a {\it term order} for $I$.

The paradigm of an ${\cal A}$-graded algebra is the semigroup algebra
$$ k [{\bf N} {\cal A}] \quad = \quad
k[ {\bf t}^{{\bf a}_1},  {\bf t}^{{\bf a}_2}, \ldots,
 {\bf t}^{{\bf a}_n}] \quad = \quad
k[{\bf x}]/I_{\cal A}.$$
The prime ideal $I_{\cal A}$ is called the
{\it toric ideal} of ${\cal A}$. It is generated by all binomials
$\,x_1^{u_1} \cdots x_n^{u_n} $ $ - x_1^{v_1} \cdots x_n^{v_n} \,$
such that $u_1 {\bf a}_1 + \cdots + u_n {\bf a}_n  =
 v_1 {\bf a}_1 + \cdots + v_n {\bf a}_n $.
We call an ${\cal A}$-graded algebra $R = k[{\bf x}]/I$  {\it coherent} if
there exists $\omega \in {\bf Z}^n$ such that $I$ is isomorphic to
 $\,in_\omega(I_{\cal A})$.

\proclaim Theorem 1.1. {\rm (Arnold, Korkina, Post and Roelofs)}
\hfill \break If $d=1$ and $n=3$ then every ${\cal A}$-graded algebra is
coherent.

Arnold [1] expressed the number of isomorphism classes of
${\cal A}$-graded algebras for ${\cal A} = \{1,p_2,p_3 \} $ in terms of the
continued fraction expansion for the
rational number $p_3/p_2$. A proof of this result
appeared in [9], and its extension to the case
${\cal A} = \{p_1,p_2,p_3\}$ was given in [10],[11].
We propose a reformulation of Theorem 1.1 using
the concept of the {\it state polytope} of a toric ideal
(see [2],[8],[16]). Recall that the
dimension of the state polytope is $n-d$. Hence for
$n=3, d=1$ it is a lattice polygon. The
edge directions of this state
 polygon  are perpendicular to the ``star'' as
defined in [9, Def.~2.9], [11, \S 2.3].

\proclaim Corollary 1.2.  If $d=1$ and $n=3$ then the isomorphism classes of
${\cal A}$-graded algebras are in bijection with the faces of the
state polytope of the toric ideal $I_{\cal A}$.

A question left open in [11] was whether these results hold for
$d=1$ and arbitrary $n$. The answer is ``no''. A first counterexample
with $n=7$ was constructed by D.~Eisenbud (personal communication).
One contribution of this paper is a new counterexample for $n=4$,
and hence a proof that Theorem 1.1 is best possible.
Incoherent monomial graded algebras for $d=1 $ and $n=4$
are the topic of Section 2 below.
In Section 3 we introduce two necessary geometric conditions for coherence,
and we construct examples of incoherent algebras which violate these
conditions. In Section 4 we characterize the radicals of
${\cal A}$-graded ideals in terms of polyhedral subdivisions.
A special case of our characterization is the main theorem
in [15] which relates triangulations and Gr\"obner bases.
In Section 5 we construct the parameter space ${\cal P}_{\cal A}$
whose points are the distinct ${\cal A}$-graded ideals in $k[{\bf x}]$.
A list of open problems and conjectures is given in Section 6.

\beginsection 2. One-dimensional mono-AGAs with four generators

An ${\cal A}$-graded algebra $R = k[{\bf x}]/I$
is called {\it monomial} (or a {\it mono-AGA}) if its
ideal $I$ is generated by monomials. The non-zero monomials
of a mono-AGA $R$ are called {\it standard}. They constitute a
$k$-vector space basis for $R$. The following is our first main result.

\proclaim Theorem 2.1. Let $d=1, n=4$ and ${\cal A}= \{1,3,4,7\}$,
and suppose $k$ is an infinite field.
\item{(a)} There exists a monomial ${\cal A}$-graded algebra
which is not coherent.
\item{(b)} There exists an infinite family of pairwise non-isomorphic
${\cal A}$-graded algebras.

\noindent {\sl Proof: }
In the  polynomial ring $k[x_1,x_2,x_3,x_4]$ we consider the
monomial ideal
$$ I \quad := \quad
 \langle \, x_1^3, \, x_1 x_2,\, x_2^2, \,x_2 x_3, \,x_1 x_4,\, x_1^2 x_3^2,
\, x_1 x_3^4,\, x_2 x_4^3,\, x_4^4 \, \rangle . \eqno (2.1) $$
The quotient algebra $\,R=k[x_1,x_2,x_3,x_4]/I \,$ is ${\cal A}$-graded.
To verify this, one must compute the Hilbert series of $I$ with respect to
the grading $\,deg(x_1)=1,\,deg(x_2)=3,\,deg(x_3)=4,\,deg(x_4)=7$.
(This can be done easily using the command {\tt hilb-numer} in the
computer algebra system MACAULAY [3].)
We list the standard monomials of low degrees:
$$
\matrix{
 1 & 2 & 3 & 4 & 5 & 6 & 7 & 8 & 9 &
10 & 11 & 12 & 13 & 14
  \cr
   x_1   &
   x_1^2   &
   x_2   &
   x_3   &
   x_1 x_3   &
   x_1^2 x_3   &
   x_4   &
   x_3^2   &
   x_1 x_3^2   &
   x_2 x_4   &
   x_3 x_4   &
   x_3^3   &
   x_1 x_3^3   &
   x_4^2
 \cr
& \cr
 15 & 16 & 17 & 18 & 19 &
20 & 21 & 22 & 23 & 24 & 25 & 26 & 27 & 28  \cr
   x_3^2 x_4   &  x_3^4   &
   x_2 x_4^2   &
   x_3 x_4^2   &
   x_3^3 x_4   &
   x_3^5   &
   x_4^3   &
   x_3^2 x_4^2   &
    x_3^4 x_4   &
   x_3^6   &
   x_3 x_4^3   &
   x_3^3 x_4^2   &
   x_3^5 x_4     &
    x_3^7 \cr} $$

\vskip .1cm

\noindent The proof is by contradiction. Suppose that $R$ is
coherent. Then there exists a rational vector
$\omega = (\omega_1,\omega_2,\omega_3,\omega_4)$
such that $I = in_\omega( I_{\cal A})$.
\item{(i)} In degree $6$ we have
$\,x_2^2 \in I \,$ but
$\, x_1^2 x_3 \not\in I$. This implies
$\, 2 \omega_2 > 2 \omega_1 + \omega_3 $.
\item{(ii)} In degree $17$ we have
$\, x_1 x_3^4 \in I \,$ but
$\, x_2 x_4^2 \not\in I$. This implies
$\,\omega_1 + 4 \omega_3 > \omega_2 + 2 \omega_4 $.
\item{(iii)} In degree $28$ we have
$\, x_4^4 \in I \,$ but
$\, x_3^7 \not\in I$. This implies
$\,4 \omega_4 > 7 \omega_3 $.

Combining these three inequalities we get
$$  (2 \omega_2) \,+\, 2 \cdot (\omega_1 + 4 \omega_3) \,+\, (4 \omega_4)
\,\,\, > \,\,\,  (2 \omega_1 + \omega_3 ) \, + \,
 2 \cdot ( \omega_2 + 2 \omega_4 ) \,+\,    ( 7 \omega_3). \eqno (2.2) $$
The left hand side and the right hand side are both equal to
$\, 2 \omega_1 + 2 \omega_2 + 8 \omega_3 + 4 \omega_4 $.
This is a contradiction, and we conclude that $R$
is not coherent. This proves part (a).

To prove part (b) of Theorem 2.1 we consider the following family of ideals:
$$ \langle \, x_1^2 x_3 - c_1 x_2^2, \, x_1 x_3^4 - c_2  x_2  x_4^2,
\,  x_3^7 - c_3  x_4^4, \,
 x_1^3 , \,
x_1 x_2,  \,
 x_1 x_4 , \,
 x_2^3  , \,
x_2^2 x_4 ,\,
 x_2 x_3,\,
 x_2 x_4^3 \,
\rangle ,
\eqno (2.3) $$
where $c_1,c_2,c_3$ are indeterminate parameters over $k$.
For every value of $c_1,c_2,c_3$ this is an ${\cal A}$-graded
ideal in $k[x_1,x_2,x_3,x_4]$. In other words, the given
three-dimensional family of ideals is flat over $k^3$.
To see this, we note that
the given generators in (2.3) are a Gr\"obner basis with respect to the
lexicographic term order induced from  $ x_1 > x_2 > x_3 > x_4 $.
Note that the three first generators in (2.3) correspond to the
three cases (i),(ii) and (iii) above. The ideal $I$ in (2.1) is
obtained from (2.3) by a deformation of the form
$\, c_1 , c_3 \rightarrow \infty , c_2 \rightarrow 0 $.

Two ideals in this family define isomorphic ${\cal A}$-graded algebras
if and only if they can be mapped into each other by an element
in the torus $(k^*)^4$  (acting naturally on the four variables).
This is the case if and only if the invariant $ c_1 c_3 / c_2^2  $
has constant value. We conclude that the ideals in (2.3) define
a one-dimensional family of non-isomorphic ${\cal A}$-graded
algebras. In particular, this family is infinite, since $k$ is infinite.
\Box

\vskip .2cm

The incoherent mono-AGA in Theorem 2.1 was
found through a  systematic search of  ${\cal A}$-graded monomial algebras.
For this search we used computational
techniques refining the ones presented in [11, \S 7].
Our point of departure was the following lemma which restricts the degrees
of minimal generators of an ${\cal A}$-graded ideal. A binomial
$\,{\bf x}^{\bf u} - {\bf x}^{\bf v} \,$
in the toric ideal $I_{\cal A}$ is called {\it primitive}
if there are no proper monomial factors
${\bf x}^{{\bf u}'}$ of ${\bf x}^{\bf u}$ and
${\bf x}^{{\bf v}'}$ of ${\bf x}^{\bf v}$ such that
$\,{\bf x}^{{\bf u}'} - {\bf x}^{{\bf v}'} \,\in \,I_{\cal A}$.
This nomenclature is consistent with [4].
Primitive binomials were called {\it star relations} in [9]
and {\it Graver basis elements} in [16].
We say that ${\bf b} \in {\bf N}{\cal A}$ is a
{\it primitive degree} if there exists
a primitive binomial of degree ${\bf b}$ in ${\cal A}$.
The following lemma for $d=1$ appears in [9, Proposition 2.1].

\proclaim Lemma 2.2. The degree of every minimal
generator of an ${\cal A}$-graded ideal  is primitive.

\noindent {\sl Proof: }
Let $I$ be an ${\cal A}$-graded ideal, let
$f$ be a  homogeneous minimal generator of $I$,
and let $\,{\bf b} := deg(f)  \in {\bf N}{\cal A}$.
We must find a primitive binomial
$\,{\bf x}^{\bf u} - {\bf x}^{\bf v}\,$
of degree ${\bf b}$ in $ I_{\cal A}$.
By the defining property (1.1),
there exists a monomial ${\bf x}^{\bf v} $
of degree ${\bf b}$ which is non-zero modulo $I$.
We may assume that $f$ has a
minimal number of monomials distinct from ${\bf x}^{\bf v}$.
Clearly, this number is at least one, that is, $f$ contains a
monomial ${\bf x}^{\bf u}$ distinct from ${\bf x}^{\bf v}$.

We claim that
${\bf x}^{\bf u} - {\bf x}^{\bf v} $ is a primitive relation in $I_{\cal A}$.
Suppose not, and let ${\bf x}^{{\bf u}'}$ be a proper  factor of
${\bf x}^{\bf u}$  and  ${\bf x}^{{\bf v}'}$ a proper
factor of ${\bf x}^{\bf v}$
such that ${\bf x}^{{\bf u}'}$ and $ {\bf x}^{{\bf v}'}$ lie in the
same graded component of $R = k[{\bf x}]/I$. Since
${\bf x}^{{\bf v}'}$ is standard, there exists $c_1 \in k$,
such that $\, {\bf x}^{{\bf u}'} -  c_1 {\bf x}^{{\bf v}'} \in I $.
By the same reasoning, there exists $c_2 \in k$
such that $\, {\bf x}^{{\bf u}-{\bf u}'} -  c_2 {\bf x}^{
{\bf v}-{\bf v}'} \in I $.
This implies $\,{\bf x}^{\bf u} - c_1 c_2 {\bf x}^{\bf v}  \in I$.
We may now replace the occurrence of ${\bf x}^{\bf u}$ in $f$
by $ c_1 c_2 {\bf x}^{\bf v} $. This is a contradiction to our minimality
assumption, and we are done. \Box

\vskip .2cm

\proclaim Proposition 2.3.
Let $d=1$ and ${\cal A} = \{a_1 < a_2 < \cdots < a_n \}
\subset {\bf N}$. Then every minimal generator of an ${\cal A}$-graded
ideal has degree at most $\, a_{n-1} \cdot a_n$.

\noindent {\sl Proof: }
It was proved in [4] that every primitive
binomial has degree at most $ a_{n-1} \cdot a_n$.
Now apply Lemma 2.2. \Box

\vskip .2cm

Proposition 2.3 improves the bound in [9, Proposition 2.11].
If $ a_{n-1}$ and $ a_n$ are relatively prime, then
the bound $\, a_{n-1} \cdot a_n \,$ is best possible.
To see this note that the binomial $\,x_n^{a_{n-1}} -
x_{n-1}^{a_n} \,$ appears in the reduced
Gr\"obner basis of $I_{\cal A}$ with respect to the
lexicographic term order induced by $x_1 \succ \cdots \succ x_n $.
The initial ideal of $I_{\cal A}$ for this term order
is an ${\cal A}$-graded ideal which has a minimal generator
of degree $\, a_{n-1} \cdot a_n $.

\vskip .1cm

The following table comprises a complete catalogue of all
incoherent mono-AGAs for $\,{\cal A} = \{a_1,a_2,a_3,
a_4 \} \,$ with $1 \leq a_1 < a_2 < a_3 < a_4 \leq 9 $.
We write the set ${\cal A}$ as a bracket $[a_1 a_2 a_3 a_4]$.
The three  integers listed immediately after each bracket are:
\item{(i)} the number of primitive binomials in $I_{\cal A}$.
\item{(ii)} the total number of all ${\cal A}$-graded monomial ideals,
\item{(iii)} the number of incoherent ${\cal A}$-graded  monomial ideals.

\noindent
If a quadruple does not appear in this list, then
all mono-AGAs are coherent for that ${\cal A}$.

\vskip .4cm

\hrule
\vskip -.1cm
$$ \matrix{
 & [1347]  &  27 &   53 &   2 &  \qquad \quad
 & [1349]  &  23 &   38 &   2 &  \qquad \quad
 & [1456]  &  26 &  51  &  2 \cr
 & [1459]  &  37 &   90 &   10  &  \qquad \quad
 & [1567]  &  35 &  79  &  6  &  \qquad \quad
 & [1568]  &  27 &   58 &   4 \cr
 & [1578]  &  33 &  79  &  2  &  \qquad \quad
 & [1678]  &  41 &  112 &  18  &  \qquad \quad
 & [1689]  &  32 &  82  &  6 \cr
 & [1789]  &  52 &  174 &   42  &  \qquad \quad
 & [2357]  &  30 &  75  &   6  &  \qquad \quad
 & [2358]  &  31 &  83  &  10 \cr
 & [2359]  &  24 &  58  &  8   &  \qquad \quad
 & [2379]  &  31 &  82  &  6  &  \qquad \quad
 & [2567]  &  30 &  67  &  2  \cr
 & [2579]  &  45 &  168 &  42 &  \qquad \quad
 & [2678]  &  27 &  53  &  2 &  \qquad \quad
 & [2689]  &  23 &  38  & 2 \cr
 & [2789]  &   41 &   113 &   10 &  \qquad \quad
 & [3459]  &  30 &  63 &  2 &  \qquad \quad
 & [3479]  &  31 &  64 &  2 \cr
 & [3578]  &  35 &  88 &  2 &  \qquad \quad
 & [3589]  &  33 &  81 &   8 &  \qquad \quad
 & [4569]  &  32 &  84 &  6 \cr
 & [4579]  &  40 &  120 &   6 &  \qquad \quad
 & [5678]  &  35 &  90  &  2 &  \qquad \quad
 & [5789]  &   40 &  113 &  2 \cr
 & [6789]  &  37 &  94  &  6 \cr} $$
\vskip -.1cm
\noindent {\bf Table 1.}
Incoherent one-dimensional mono-AGAs
with $n=4$ and degrees $\leq 9 $.

\vskip ,2cm

\hrule

\vskip .4cm

\beginsection 3. Polyhedral conditions for coherence

The computational results in Section 2 raise the question
whether there exist structural features
of coherent AGAs which are not shared by all AGAs.
Here we  identify two such features:
standard monomials (Theorem 3.3)
and degrees of minimal generators (Theorem 3.6)
are subject to certain geometric  restrictions
in the coherent case.
Our presentation  assumes familiarity with the language
of combinatorial convexity (see e.g.~[6],[17]).

For a fixed set ${\cal A} = \{{\bf a}_1,\ldots,{\bf a}_n\}
\subset {\bf N}^d \! \setminus \! \{0\}$ we consider the linear map
$$ deg \, : \,{\bf N}^n \,\rightarrow \,{\bf N} {\cal A} \,, \,\,\,
{\bf u} = (u_1,\ldots,u_n) \,\mapsto \,\sum_{i=1}^n u_i \cdot {\bf a}_i.$$
Each inverse image $\,deg^{-1}({\bf b}) \,$ consists of (the
exponent vectors of) the monomials of degree ${\bf b}$.
Following [16] we form their convex hull
$\,P[{\bf b}] := conv( deg^{-1}({\bf b}))$.
This polytope is called the {\it fiber} of ${\cal A}$ over
${\bf b}$. By the fiber of a monomial
$\,{\bf x}^{\bf u}\,$ we mean $\,P[ deg({\bf u})]$.

\proclaim Observation 3.1.
Every standard monomial ${\bf x}^{\bf u}$ of a coherent mono-AGA
corresponds to a vertex ${\bf u}$ of its fiber $\,P[deg({\bf u})]$.

\noindent {\sl Proof:}
Let ${\bf x}^u$ be standard in $k[{\bf x}]/in_\omega (I_{\cal A})$.
Then ${\bf u}$ is the unique point
in $ deg^{-1}(deg({\bf u}))$ at which the linear functional
$\omega$ attains its minimum. Hence  ${\bf u}$ is a vertex
of $P[deg({\bf u})]$. \Box

\vskip .2cm

The theory of ${\cal A}$-graded algebras provides an abstract
setting for the study of integer programming problems with
respect to a fixed matrix (cf.~[16]).
This was a main motivation for writing the present paper.
The subsequent remark makes it more precise.

\proclaim Polyhedral Remark 3.2. {\rm
Consider the following alternative definition for the object of study
in Section 2: A {\it mono-AGA} is a rule which selects
one lattice point (called standard monomial) from each fiber,
subject to the axiom that the set of
standard monomials is closed under divisibility.
This selection is {\it coherent} if it is induced by a linear
functional $\omega$. Our choice of the term ``coherent''
parallels the notions of {\it coherent subdivisions} and
{\it coherent triangulations}. We refer to [6],[17] and the references
given there (and to Section 3 below). Also our usage of the letter
``${\cal A}$''  is consistent with that of ~[6], [15].}

\vskip .1cm

Does there exist an incoherent mono-AGA which has a standard
monomial that is not a vertex of its own fiber~?
The answer was found to be ``no'' for all $218$ incoherent mono-AGAs
listed in Table~1. We do not know the answer for
$d=1$ and $n=4$ in general. For $d=1$ and $n=5$ we can show that
the answer is ``yes''.

\proclaim Theorem 3.3. Let $d \! = \! 1, n
\! = \! 5$ and ${\cal A}= \{ 3,4,5,13,14\}$.
There exists a monomial ${\cal A}$-graded algebra with a
standard monomial ${\bf x}^{\bf u}$ such that
${\bf u}$ is not a vertex of its fiber $P[deg({\bf u})]$.

\noindent {\sl Proof. }  In the polynomial ring
$k[x_1,x_2,x_3,x_4,x_5]$ we consider the ideal
$$ I \quad = \quad
\langle \,x_1^3,\, x_2^2, \, x_3^2,\, x_1 x_5,\, x_2 x_5,
   \,x_3 x_5, \, x_5^2 \, \rangle . $$
This ideal is ${\cal A}$-graded. Indeed,
an easy MAPLE or MACAULAY computation shows that
$R = k[{\bf x}]/I$ has the correct Hilbert series, namely,
$$ {1 \over 1 - t} - t - t^2 \quad = \quad
\sum_{m \in {\bf N} {\cal A}} t^m . $$
The monomial $\,x_1^2 x_2 x_3 \,$ does not lie in $I$: it is a
standard monomial of degree $15$. Its fiber is the Newton polytope
of the coefficient $\,c_{15} ({\bf x})\,$ in the expansion of
the formal power series
$$  1  / \bigl( (1 - x_1 t^3) (1- x_2 t^4) (1- x_3 t^5)
(1- x_4 t^{13}) (1- x_5 t^{14}) \bigr)
\quad = \quad
\sum_{m=0}^\infty \, c_m (x_1,x_2,x_3,x_4,x_5) \cdot t^m .$$
We find that $\,c_{15} ({\bf x}) \, = \,
x_1^5 + x_1 x_2^3 + x_1^2 x_2 x_3 + x_3^3 $.
The point $\,(2,1,1,0,0)\,$ lies in the relative interior
of the triangle
$\,P[(2,1,1,0,0)] \,= \, conv \,\{
(5,0,0,0,0), (1,3,0,0,0), (0,0,3,0,0) \} $. \Box

\vskip .2cm

After this discussion of vertices of fibers, we now turn
our attention to edges of fibers.
An element ${\bf b}$ of ${\bf N}{\cal A}$
is said to be a {\it Gr\"obner degree}
if a binomial  of degree ${\bf b}$ appears in some
reduced Gr\"obner basis of $I_{\cal A}$.
Equivalently, the Gr\"obner degrees are precisely the
degrees of the minimal generators of all {\sl coherent}
${\cal A}$-graded ideals.
Lemma 2.2 implies that every Gr\"obner degree is also a primitive degree.
The following geometric characterization of
Gr\"obner degrees was proved in the joint work  [16]
with R.~Thomas.

\proclaim Theorem 3.4.
An element ${\bf b}$ of ${\bf N}{\cal A}$
is a {\it Gr\"obner degree} if and only if its fiber
$P[{\bf b}] $ has an edge which is not parallel to
any edge of a different fiber $P[{\bf b'}]$ with
${\bf b'} \leq {\bf b}$.

\noindent {\sl Proof: }
This is a reformulation of Theorem 5.1 in [16]. \Box

\vskip .1cm

Theorem 3.4 implies that a primitive binomial lies
in some reduced Gr\"obner basis of $I_{\cal A}$
if and only if its two monomials
are connected by an edge in their fiber.

It was shown in [16, \S 5] that there may exist
primitive degrees which are not Gr\"obner.
Here is an example for $d=1$. It is derived from
Korkina's example in [9, Remark 3.3].

\vskip .3cm

\noindent {\bf Example 3.5.}
{\sl (A primitive degree which is not a Gr\"obner degree)  }
\hfill \break
 Let ${\cal A} =  \{ 15,20,23,24 \}$.
The binomial $\,x_1^2 x_2^3 x_4^2 - x_3^6 \,$ is
the unique primitive binomial of degree $\,138 $.
(This can be verified easily using [16, Algorithm 4.3].) \
We shall prove
that it is not an element of any reduced Gr\"obner basis of $I_{\cal A}$.
Consider the convex combination
$$  {1 \over 4} (0,0,6,0) \,+ \,
    {3 \over 4} (2,3,0,2) \quad = \quad
 {1 \over 4} (5,2,1,0)  \,+ \,
{1 \over 4}  (1,5,1,0) \,+ \,
 {1 \over 2} (0,1,2,3) .$$
This shows that $\,conv \{(2,3,0,2), (0,0,6,0) \}$,
the segment corresponding to our binomial, is not an
edge of its fiber $P[138]$. Using Theorem 3.4
we conclude that the given binomial does not lie
in any reduced Gr\"obner basis.
Therefore $138$ is a primitive degree which is not
a Gr\"obner degree. \Box

\vskip .2cm

This raises the question whether  there exists some
(necessarily incoherent) ${\cal A}$-graded algebra
which has a minimal generator of non-Gr\"obner degree.
The answer is ``no'' for the set ${\cal A} =  \{ 15,20,23,24 \}$
above, but it becomes ``yes'' after adding sufficiently
many new generators. Clearly, the number $n=145$ below
is not best possible.

\proclaim Theorem 3.6.
For $d = 1$ and $n = 145 $
there exists an ${\cal A}$-graded ideal $I$ which has
a minimal generator whose degree is not a Gr\"obner degree.

\noindent {\sl Proof: }
Let ${\cal A}' = \{15,20,23,24,107,109 \}$.
Let $S'$ be a polynomial ring in six variables
$x_{15}, x_{20}, x_{23}, x_{24}, x_{107}, x_{109}$.
We grade $S$ by setting $\,deg(x_i) = i $.
Let $M'$ be the ideal generated by the six variables, and let
$M'_{\geq 139}$ be the ideal generated by all
monomials of degree $\geq 139$ in $S'$.
Let $I$ be the binomial ideal generated by
$$ \eqalign{
& x_{15}^4 ,\,\,
x_{20}^2 x_{23},\,\,
x_{15}^3 x_{23},\,\,
x_{15}^3 x_{24},\,\,
x_{24} x_{23}^2,\,\,
x_{15} x_{23}^3,\,\,
x_{20} x_{24}^3,\, \,
 x_{15} x_{20}^4,\,\,
x_{15}^2 x_{20} x_{23}^2,\,  \cr &
x_{24}^5,\,\,
x_{15} x_{93},\,\,
x_{15} x_{109},\,\,
x_{24} x_{109},\,\,
x_{15}^2 x_{107},\,\,\,\hbox{and} \,\, \,
x_{15}^2 x_{20}^3 x_{24}^2-x_{23}^6. \cr}
$$
The binomial $x_{15}^2 x_{20}^3 x_{24}^2 - x_{23}^6$ has degree
$138$, and this degree is primitive but not Gr\"obner (by Example 3.5).
The ideal $I$ is constructed to have the following property:
the Artinian ring  $S'/(I + M'_{\geq 139})$ is ${\cal A}$-graded
up to degree $138$. In other words, its Hilbert series equals
$\, \sum \bigl\{\,t^b \,:\, b \in {\bf N} {\cal A} \,\,\hbox{and}\,\,
b \leq 138 \bigr\} $.

Let $\,{\cal A}'' \, = \, \{139,140,\ldots,277 \}\,$ and introduce
the corresponding polynomial ring
$\,S'' = k [x_{139}, x_{140}, \ldots,x_{277}]$.
We write  $M''$ for the ideal generated by all $139$ variables in $S''$,
and we let $J''$ be any ${\cal A}''$-graded ideal in $S''$.

Finally, we set ${\cal A} := {\cal A}' \cup {\cal A}''$, and
we introduce the corresponding $145$-variate polynomial ring $\, S \, :=
 \, S' \otimes_k S'' $. In this ring  we form the ideal
$$ J \quad := \quad
\langle M' \cdot M'' \rangle \,+\,
\langle  I + M'_{\geq 139}\rangle \,+\,
\langle J'' \rangle. $$
By construction, the ideal $J$ is ${\cal A}$-graded in $S$.
It has the primitive binomial
$\,x_{15}^2 x_{20}^3 x_{24}^2-x_{23}^6 \,$ among its
minimal generators. However, its degree $138$ is not
a Gr\"obner degree for ${\cal A}$. This completes the proof.
\Box

\vskip .3cm

\beginsection 4. The radical of an ${\cal A}$-graded ideal

In this section	 a polyhedral construction is presented
for the radical of  any ${\cal A}$-graded ideal.
The {\it positive hull} of ${\cal A}$ is the closed convex
polyhedral cone
$$ \,pos({\cal A}) \quad := \quad
\bigl\{ \,\sum_{i=1}^n \lambda_i \cdot {\bf a}_i\,
 \,:\,(\lambda_1,\ldots,\lambda_n) \geq 0 \,\bigr\}. $$
If $\sigma$ is any subset of ${\cal A}$ then we similarly write
$pos(\sigma)$ for the positive hull of $\sigma$.
By a {\it face} of $\sigma$ we mean a subset $\tau$
such that the cone $\,pos( \tau ) \,$
is a face of the cone $\,pos(\sigma) \,$
and $\tau = pos ( \tau) \cap \sigma $. We identify the toric
ideal $I_\sigma$ corresponding to a subset $\sigma$ with the prime
ideal $\,I_\sigma + \langle x_i : a_i \not\in \sigma \rangle\,$
in  $k [{\bf x}]$.
A {\it polyhedral subdivision} of ${\cal A}$ is a
collection $\Delta $ of subsets of ${\cal A}$ such that
$\,\{ \, pos(\sigma) \,:\, \sigma \in \Delta \}\,$ is a
polyhedral fan with support equal to the cone $ pos ({\cal A})$.

A basic construction due to R.~Stanley associates
to any integral polyhedral complex a radical binomial ideal.
(See [5, Example 4.7] for an algebraic discussion.)
If $\Delta$ is a polyhedral subdivision of ${\cal A}$, then
its {\it Stanley ideal} is $\,I_\Delta \,:= \,
\bigcap_{\sigma \in \Delta} I_\sigma$. We remark
that $ R = k[{\bf x}]/I_\Delta $ is also graded by the semigroup
${\bf N} {\cal A}$, but it is generally \underbar{not}
${\cal A}$-graded because, for some ${\bf b} \in {\bf N} {\cal A}$,
the graded component $R_{\bf b}$ may be zero.
 Finally, we call two arbitrary ideals $I$ and $I'$ in  $k [{\bf x}]$
 {\it torus isomorphic} if there exists
$\lambda \in (k^*)^n$ such that (1.2) holds.
The following is the main result in this section.

\proclaim Theorem 4.1.
If $I$ is any ${\cal A}$-graded ideal,  then there
exists a polyhedral subdivision $\Delta$ of ${\cal A}$
such that $\,Rad(I) = \bigcap_{\sigma \in \Delta} J_\sigma \,$
where each ideal $J_\sigma$ is prime and  torus isomorphic to $I_\sigma$.

We shall make a two remarks before presenting the proof.
First, as a special case of Theorem 4.1, we can recover the
main result in [15]. Namely, this is the case where
\item{$\bullet$} the given ideal $I$ is a coherent mono-AGA, and
\item{$\bullet$} the semigroup ${\bf N}{\bf A}$ is graded
(or, equivalently, the toric ideal $I_{\cal A}$ is homogeneous).

\noindent If these two hypotheses are met, then
$pos({\cal A})$ is the cone over
a polytope $conv({\cal A})$, and  $\Delta$  is (the
complex of cones over) a  {\it regular  triangulation} of this polytope.

Our second remark is to explain the mysterious appearance of
the ideals $J_\sigma$. What is
the point of replacing $I_\sigma$ by a torus isomorphic ideal $J_\sigma$,
for each maximal cell $\sigma$ of $\Delta$~?
The answer is that the Stanley ideal
$\,I_\Delta = \bigcap_{\sigma \in \Delta} I_\sigma \,$ itself
is generally \underbar{not} torus isomorphic to the radical of $I$.
We present an example where this happens.

\vskip .2cm

\noindent {\bf Example 4.2.}
{\sl An ${\cal A}$-graded ideal $I$ such that
$Rad(I)$ is not torus isomorphic to any $I_\Delta$. } \break
\rm
Let $d=3, n=6$ and
${\cal A} =
\{(4,0,0), (0,4,0),(0,0,4),   (2,1,1),(1,2,1),(1,1,2) \}$.
For every choice of non-zero constants $c_1,c_2,c_3 \in k^*$,
the following ideal is ${\cal A}$-graded:
$$ \eqalign{
 I_{c_1,c_2,c_3} \quad = \quad
   \langle & x_1 x_2 x_3,\, x_1 x_5 x_6,\, x_2 x_4 x_6,
\, x_3 x_4 x_5,\,    x_1 x_2 x_6^2,\, x_1 x_3 x_5^2,\, x_2 x_3 x_4^2,\, \cr
   & \,\,\,\, x_1 x_5^4 - c_1 x_2 x_4^4 ,\,\,
    x_2 x_6^4 - c_2 x_3 x_5^4, \,\,
    x_3 x_4^4 - c_3 x_1 x_6^4 \, \rangle. \cr }$$
The underlying subdivision $\Delta$ of ${\cal A}$
consists of three quadrangular cones and one
triangular cone. This can
be seen from the prime decomposition
$\, Rad( I_{c_1,c_2,c_3} ) \, = $
$$ \langle x_1,x_2,x_3 \rangle \,\cap \,
   \langle x_1 x_5^4 - c_1 x_2 x_4^4 ,x_3,x_6 \rangle \,\cap \,
   \langle x_2 x_6^4 - c_2 x_3 x_5^4 ,x_1,x_4 \rangle \,\cap \,
   \langle x_3 x_4^4 - c_3 x_1 x_6^4 ,x_2,x_5 \rangle . $$
{}From this decomposition we can see that $Rad(I_{c_1,c_2,c_3})$ is
torus isomorphic to the Stanley ideal $I_\Delta$ if and only if
the invariant $\,c_1 c_2 c_3 \,$ attains the value $1$.
The reader will not fail to note a certain analogy between this example
and the three-dimensional family in (2.3).

The reason for the phenomenon in Example 4.2 is the existence
of moduli (infinite families) of isomorphism classes of
${\cal A}$-graded algebras.  These moduli stem from
``extraneous components'' in the parameter space ${\cal P}_{\cal A}$
(see Section 5 and Problem 6.4).

\vskip .1cm

For the proof of Theorem 4.1 we shall need the following lemma.

\proclaim Lemma 4.3.  Let $I \subset k[{\bf x}]$ be an
${\cal A}$-graded ideal which contains no monomials.
Then $I$ is isomorphic to the toric ideal $I_{\cal A}$.

\noindent {\sl Proof: } This follows from the
characterization of Laurent binomial ideals in [5, \S 2]. \Box

\vskip .3cm

\noindent {\sl Proof of Theorem 4.1: }
The polynomial ring $S = k[x_1,\ldots,x_n]$
is graded by the semigroup ${\bf N} {\cal A}$
via $\,deg(x_i) = {\bf a}_i $. For any ${\bf b} \in {\bf N} {\cal A}$
we define the subalgebra  $\,S_{({\bf b})} := \bigoplus_{m=0}^\infty
S_{m {\bf b} }$. This algebra is generated by a finite set of monomials.
Inside it we consider the binomial ideal $\, I_{({\bf b})} := I \cap
S_{({\bf b})}$. The corresponding subalgebra $R_{({\bf b})}
:= S_{({\bf b})}/I_{({\bf b})}$ of our given
${\cal A}$-graded algebra $R=S/I$ is a finitely
generated $k$-algebra of Krull dimension $1$.
It is not possible that all elements in such an algebra are nilpotent.
We conclude that there exists a monomial ${\bf x}^{\bf u}$ in
$S_{({\bf b})}$ which is not nilpotent modulo $I_{({\bf b})}$.

Let ${\bf x}^{\bf u}$ and ${\bf x}^{\bf v}$ be two such
non-nilpotent monomials in $R_{({\bf b})}$.
We claim that their product ${\bf x}^{\bf u} {\bf x}^{\bf v} \in
R_{({\bf b})}$ is not nilpotent either. To see this,
we choose integers $m_1$ and $m_2$ such that
$ x^{m_1 {\bf u}}$ and $x^{m_2 {\bf v} }$ have the same degree.
There exists a non-zero constant $c \in k^*$ such that
$ \,{\bf x}^{m_1 {\bf u}} = c \cdot {\bf x}^{m_2 {\bf v}} $
in $R_{({\bf b})}$. This implies
$\,({\bf x}^{m_1 {\bf u}} {\bf x}^{m_2 {\bf v}})^m = c^m ({\bf x}^{
{\bf v}})^{2 m m_2 } = c^{-m} ({\bf x}^{\bf u})^{2m m_1} \not= 0 \,$
in $R$ for all $m > 0 $,  and consequently $({\bf x}^{\bf u}
{\bf x}^{\bf v})^m \not=
0 $ in $R_{({\bf b})}$ for all $m > 0 $.
We have shown that the set of non-nilpotent monomials in
$R_{({\bf b})}$ is multiplicatively closed.

This multiplicativity property allows us to synthesize
the polyhedral subdivision $\Delta$. The {\it support} of a
monomial is defined as  $\, supp({{\bf x}^{\bf u}}) \,:= \,
\{ \,{\bf a}_i \in {\cal A}\,:\,u_i \not= 0 \}$.
Clearly, we have $\,supp({\bf x}^{\bf u} {\bf x}^{\bf v})
\,=\, supp({\bf x}^{\bf u}) \, \cup \, supp({\bf x}^{\bf v}) $.
This implies that the set of supports of
non-nilpotent monomials in $ R_{({\bf b})}$ has a unique
maximal element. This subset of ${\cal A}$ is denoted by
$\, cell({\bf b})$.  We define $\Delta$ to be the
collection of all subsets  $cell({\bf b})$ as
${\bf b}$ ranges over ${\bf N} {\cal A}$.

We shall prove that $\Delta$ is indeed a polyhedral subdivision
of ${\cal A}$. Let $\tau$ be any face of
$\sigma = cell({\bf b})$  (possibly $\tau = \sigma$),
and let ${\bf b}'$ be any lattice point in the
relative interior of $pos(\tau)$.
It suffices to show that $\,cell({\bf b}') = \tau $.
By the property of being a face, $\tau$ is the unique
maximal subset of $\sigma$ which is the support of
any monomial ${\bf x}^{{\bf u}'}$ in  $S_{({\bf b}')}$.
Such a monomial $ {\bf x}^{{\bf u}'}$ is not nilpotent modulo $I$,
since there exists a monomial ${\bf x}^{\bf u}$ of degree ${\bf b}$
whose support equals $\sigma \supset \tau =
supp({\bf x}^{{\bf u}'})$.  Suppose there exists a
non-nilpotent monomial ${\bf x}^{{\bf u}''}$ in $R_{({\bf b}')}$
whose support $\, \rho := supp( {\bf x}^{{\bf u}''})\,$ properly
contains $\tau $. Then $\,\rho \setminus \sigma = \rho
\setminus \tau \not= \emptyset $. Choose integers $m_1, m_2$
and a non-zero constant $c \in k^*$ such that
${\bf x}^{m_1 {\bf u}'} - c \cdot {\bf x}^{m_2 {\bf u}''} \in I $.
Let the degree of this binomial be $ m_3 \cdot {\bf b}' $.
Choose an integer $m_4 \gg 0 $ such that $\, m_4 \cdot {\bf b}
- {\bf b}' \,$ lies in the relative interior of
$pos( \sigma)$, and let ${\bf x}^{\bf w}$ be a monomial
having degree $m_4 \cdot {\bf b} - {\bf b}'$ and support $\sigma$.
We conclude that  $ x^{ m_3 {\bf w} + m_1 {\bf u}'}
- c \cdot x^{m_3 {\bf w} + m_2 {\bf u}''}\,$ lies in the degree
$m_3 m_4 \cdot {\bf b} $ component of the ideal $I$.
The first monomial $ {\bf x}^{ m_3 {\bf w} + m_1 {\bf u}' }$
is not nilpotent modulo $I$ since it has support $\sigma$.
The second monomial $ {\bf x}^{m_3 {\bf w} + m_2 {\bf u}''}$
is nilpotent modulo $I$ since its support
$\sigma \cup \rho$ strictly contains $\sigma$.
This is a  contradiction, and we conclude
that $cell({\bf b}') = \tau$.
This completes the proof that $\Delta$ is
a polyhedral subdivision of ${\cal A}$.

We next compute the radical of $I$.
Let $\sigma $ be a maximal cell in $\Delta$.
By construction, the elimination ideal
$\,I \cap k[ x_i : i \in \sigma] \,$ is
a $\sigma$-graded ideal which contains no monomials.
Lemma 4.3 implies that its natural embedding into $k[{\bf x}]$,
$$\,J_\sigma \,\,\,\,:=\,\,\,\, ( I \cap k[ x_i : i \in
\sigma]) \,+ \, \langle x_j \,:\, a_j \not\in \sigma \rangle, $$
is torus isomorphic to the toric prime $I_\sigma$.
We claim that $ \,Rad(I)\, = \,
\bigcap_{\sigma \in \Delta} J_\sigma $.

We first show the inclusion
$\, I \subseteq \bigcap_{\sigma \in \Delta} J_\sigma $.
(This automatically implies
$\, Rad(I) \subseteq \bigcap_{\sigma \in \Delta} \! J_\sigma \,$
because the right hand side is radical.)
If ${\bf x}^{\bf u}$ is any monomial not contained in
$\bigcap_{\sigma \in \Delta} \langle x_j : j \not\in \sigma \rangle $,
then $\,supp({\bf x}^{\bf u}) \subseteq \sigma$ for some $\sigma
\in \Delta$, and hence $ {\bf x}^{\bf u}$ is not nilpotent
modulo $I$. This shows that all monomials which are nilpotent
modulo $I$ lie in $\bigcap_{\sigma \in \Delta} J_\sigma $. Consider
any binomial $\,f := {\bf x}^{\bf u} - c \cdot {\bf x}^{\bf v} $ in $I$
with  both terms not nilpotent modulo $I$.
Let ${\bf b} = deg({\bf x}^{\bf u}) = deg({\bf x}^{\bf v})$.
Fix $\sigma \in \Delta$. If $cell({\bf b})$ is
a face of $\sigma$, then $\,f \,\in \,I_{({\bf b})} \cap
k[ x_i : i \in \sigma] \,\subset\, J_\sigma$. If $cell({\bf b})$
is not a face of $\sigma$, then  the supports of ${\bf x}^{\bf u}$
and ${\bf x}^{\bf v}$ are not subsets of $\sigma$.
Therefore both ${\bf x}^{\bf u}$ and ${\bf x}^{\bf v}$ lie in
$ \langle x_j : j \not\in \sigma \rangle $, and hence $f \in J_\sigma$.

For the reverse inclusion $ \,\bigcap_{\sigma \in \Delta}
J_\sigma  \subseteq Rad(I) \,$ we use the Nullstellensatz:
it suffices to prove that the variety $\,{\cal V}(I)\,$
is contained in $\,\cup_\sigma {\cal V}(J_\sigma)$.
Let ${\bf u} \in {\bar k}^n $ be any zero of $I$,
where $\bar k $ is the algebraic closure of $k$.
Abbreviate  $\,\rho :=  supp({\bf u})$.
Consider the monomial $\prod_{i \in \rho} x_i $
and let ${\bf b}$ be its degree. Let $\sigma$ be any
maximal cell of $\Delta$ which has $cell({\bf b})$ as a face.
By construction, no power of  $\prod_{i \in \rho} x_i $
vanishes at ${\bf u}$. Hence the monomial
$\prod_{i \in \rho} x_i $ is not nilpotent modulo $I$, and its
support $\rho $ is a subset of $cell({\bf b}) \subseteq \sigma$.
In other words, ${\bf u}$ is a zero of the ideal
$\,\langle x_j : j \not\in \sigma \rangle $. Therefore
${\bf u}$ is a zero of $J_\sigma$.
This completes the proof. \Box

\vskip .3cm
\beginsection 5. The parameter space.

In this section we aim to answer the question
posed in the first sentence of the introduction.
We construct the parameter space ${\cal P}_{\cal A}$
whose points are in bijection with the distinct
${\cal A}$-graded ideals in $k[{\bf x}]$. The torus $(k^*)^n$ acts
naturally on the space ${\cal P}_{\cal A}$, and its orbits are in
bijection with the isomorphism types of ${\cal A}$-graded algebras.

Let ${\cal A} = \{{\bf a}_1,\ldots,{\bf a}_n\} \subset {\bf N}^d
\setminus \{ 0 \}$ as before,
and fix the ${\bf N}{\cal A}$-grading of $S = k[{\bf x}]$
given by $deg(x_i) = {\bf a}_i$.
For any integer $r > 0$ we define the following {\it zonotope} in ${\bf R}^d$:
$$ Z_r ({\cal A}) \quad := \quad \bigl\{
\,\sum_{i=1}^n \lambda_i \cdot {\bf a}_i \,\,:\,\,
0 \leq \lambda_i \leq r \,\,\hbox{for}\,\,i=1,\ldots,n \,\bigr\}.
\eqno (5.1) $$
Let $M_r$ be the ideal in $S = k[{\bf x}]$ spanned by all monomials
${\bf x}^{\bf u}$ such that $\,deg({\bf x}^{{\bf u} +{\bf v}}) \not\in
{\cal Z}_r ({\cal A})\,$ for all ${\bf v} \in {\bf N}^n $.
We write $S^{(r)}$ for the quotient algebra $S/M_r$.
The algebra $S^{(r)}$ is graded by ${\bf N} {\cal A}$, and it is
artinian because each variable $x_i$ has a power lying in $M_r$.
An ideal $J \subset S^{(r)}$ is called {\it ${\cal A}$-graded} if
$\,dim_k \bigl( (S^{(r)} /J )_{\bf b} \bigr) \, = \, 1\,$ for all
${\bf b} \in Z_r ({\cal A}) \cap {\bf N}{\cal A}$.
If $I$ is any ${\cal A}$-graded ideal in $S$, then its image
$\,I^{(r)} \, := \, (I + M_r) /M_r\,$ is an ${\cal A}$-graded ideal
in $S^{(r)}$.

\proclaim Proposition 5.1. There exists an integer $r \gg 0$ such that the
assignment
$\, I \,\mapsto \, I^{(r)}\,$ defines a bijection between the ${\cal A}$-graded
ideals in $S$ and the ${\cal A}$-graded ideals in $S^{(r)}$.

\noindent {\sl Proof: } For the injectivity of the  map
$\,I \,\mapsto \, I^{(r)} \,$ it suffices to choose $r$ such that the
zonotope $Z_r({\cal A})$ contains all primitive degrees.
For instance, if $a$ is the maximum of the Euclidean norms $||{\bf a}_i ||$,
then $\,r = (n-d) \cdot  a^d \,$ has this property by [15, \S2].
If  $ I $ and $J$ are two distinct ${\cal A}$-graded ideals
in $S$, then, by Lemma 2.2, there exists a primitive degree ${\bf b}$
such that $I_{\bf b}$ and $J_{\bf b}$ are distinct  subspaces of
$S_{\bf b}$. But $S_{\bf b} = S^{(r)}_{\bf b}$, hence
$\,I^{(r)}_{\bf b} = I_{\bf b} \,\not=\, J_{\bf b} = J^{(r)}_{\bf b} $,
and therefore $I^{(r)}$ and $J^{(r)}$ are distinct ${\cal A}$-graded
ideals in $S^{(r)}$.

To prove surjectivity we choose $r \gg 0$ to have the following property:
If $I$ is any binomial ideal in $k[{\bf x}]$ whose generators
have primitive degrees, and $\prec$ is any term order, then
the initial ideal $in_\prec(I)$ is generated by monomials of degree $r$.
For instance, combining the above bound with the
known doubly-exponential degree bounds for Gr\"obner bases [14], we
see that  the choice $\, r \,=\,  (n-d)^{2^n} \cdot a^{d 2^n}\,$
will surely suffice.

Let $\,J \, = \,\langle f_1,\ldots,f_m \rangle\,$ be any
${\cal A}$-graded ideal in $S^{(r)}$. Here the $f_i$ are binomials
of primitive degree, so they have a unique preimage in $S$. We consider
the ideal in $S$ generated by these preimages
$\,I\,:=\,\langle f_1,\ldots,f_m \rangle $. Since the
conclusion $I^{(r)} = J$ is automatic, we only have to show that
$I$ is ${\cal A}$-graded. Equivalently, we must show that $I$ has the
same ${\bf N}{\cal A}$-graded Hilbert series as the toric ideal
$I_{\cal A}$. Let $\prec$ be any term order. Since the Hilbert series
is preserved under passing to the initial monomial ideal, it suffices
to show that $in_\prec(I_{\cal A})$ and
$in_\prec(I)$ have the same Hilbert series.
Let $\, {\bf x}^{{\bf u}_1} ,\ldots,{\bf x}^{{\bf u}_{s}} \,$
be the minimal generators of $in_\prec(I)$.
The  Hilbert series of interest equals the rational function
$$ H ( I; {\bf t} ) \,\, =\,\,
H( in_\prec(I) ; {\bf t}) \quad = \quad
{ \sum_{ \nu \subseteq \{1,\ldots,s\}} (-1)^{\vert \nu \vert} \cdot
deg \bigl( lcm (\{ {\bf x}^{{\bf u}_j} \,:\, j \in \nu \})\bigr)
\over \prod_{i=1}^n  (1 - t_1^{a_{1i}} \cdots t_d^{a_{di}} )}.
\eqno (5.2) $$
By construction, the degree of each term in the numerator polynomial
lies in $Z_r({\cal A})$, and the same is true for $I_{\cal A}$.
Therefore $\,H(I,{\bf t}) - H(I_{\cal A},{\bf t})\,$ is a rational
function of the form
$\,p({\bf t})/\prod_{i=1}^n  ( 1 - {\bf t}^{{\bf a}_i})$, where
all monomials appearing in $p({\bf t})$ lie in $ Z_r ({\cal A})$.
Moreover, since $I$ is ${\cal A}$-graded in $S^{(r)}$, no multiple
of a monomial appearing in
the power series expansion of $\,H(I,{\bf t}) - H(I_{\cal A},{\bf t})\,$
can lie in $Z_r({\cal A})$.  In other words, the image of
$\,H(I,{\bf t}) - H(I_{\cal A},{\bf t})\,$ in the
artinian ring $\,k[t_1,\ldots,t_d]/
\langle \,{\bf t}^{\bf b} \,:\,{\bf b} \not\in Z_r ({\cal A}) \,\rangle \,$
is zero.  Indeed, in this ring the product
$\prod_{i=1}^n  ( 1 - {\bf t}^{{\bf a}_i})$ is invertible, and
we can conclude that $p({\bf t})$ is the zero polynomial.
\Box

\vskip .2cm

Proposition 5.1 is somewhat unsatisfactory in that the
doubly-exponential lower bound for $r$
used in its proof seems too big. We conjecture that
the choice $\,r = (n-d) \cdot  a^d \,$ is large enough.
In fact, our argument shows that this choice is large enough to give
the desired bijection for ${\cal A}$-graded {\sl monomial} ideals.
However, even for monomial ideals, it is not enough to require
``${\cal A}$-gradedness'' only up to the primitive degrees.

\vskip .2cm

\noindent {\bf Example 5.2.} {\sl
(One-dimensional in all primitive degrees does not imply
${\cal A}$-graded) } \break
Let ${\cal A} = \{(3,0),(2,1),(1,2),(0,3) \}$. The toric ideal
$I_{\cal A}$ is ideal of the twisted cubic curve in $P^3$.
It contains precisely five primitive binomials \ (cf.~[15, \S 4]):
$$ x_1 x_3 - x_2^2, \,\, x_1 x_4 - x_2 x_3, \,
\, x_2 x_4 - x_3^2 ,\,\, x_1^2 x_4 - x_2^3,\,
\hbox{ and } \, x_1 x_4^2 - x_3^3 .$$
The set of primitive degrees is
$\, {\cal D} \, = \, \{ (4,2),(3,3),(2,4), (6,3),(3,6) \}$.
Consider the monomial ideal
$\,\, I \, := \, \langle \, x_1 x_4,\, x_2^2,\, x_3^2 \, \rangle $.
The quotient algebra $\,k[{\bf x}]/I \,$ is one-dimensional in all degrees
${\bf b}$ with ${\bf b} \leq {\bf d}$ for some ${\bf d} \in {\cal D}$.
But $I$ is not
${\cal A}$-graded since it contains \underbar{all} monomials
of degrees $(4,5)$ and $(5,4)$. \Box

\vskip .2cm

We are now prepared to construct the parameter space of
${\cal A}$-graded ideals. Let ${\bf b} \in {\bf N}{\cal A}$.
Consider the vector space $\,k^{deg^{-1}({\bf b})}\,$ of all $k$-valued
functions on the fiber $\,deg^{-1}({\bf b})$, and let
$\,{\bf P}_{\bf b}\,$ denote the projectivization of
$\,k^{deg^{-1}({\bf b})}$.
Choose $r \gg 0$ as in Proposition 5.1 and form the product
of projective spaces $\,{\bf P} \,:= \,\Pi_{\bf b} {\bf P}_{\bf b}$, where
${\bf b}$ runs over all points in $Z_r({\cal A}) \,\cap\, {\bf N}{\cal A}$.
If $f$ is an element of the product ${\bf P}$, then
we write $f^{\bf b} \in {\bf P}_{\bf b}$ for its ${\bf b}$-th component,
and, if  ${\bf u} \in {\bf N}^n$ with $deg({\bf u}) = {\bf b}$, then
$\,f^{\bf b}_{\bf u} \,$
denotes the homogeneous coordinate of $f^{\bf b}$ indexed by ${\bf u}$.
We define a closed subscheme ${\cal P}_{\cal A}$ of ${\bf P}$ by
the equations
$$ \quad f^{\bf b}_{\bf u} \cdot f^{{\bf b} + {\bf c}}_{{\bf v} + {\bf w}}
\quad = \quad
   f^{\bf b}_{\bf v} \cdot f^{{\bf b} + {\bf c}}_{{\bf u} + {\bf w}}\,
\qquad  \hbox{whenever} \,\,\,\,
deg({\bf u}) = deg({\bf v}) = {\bf b} \,\,
\hbox{and} \,\, deg({\bf w}) = {\bf c}. \eqno (5.3) $$
The following is the main result in this section.

\proclaim Theorem 5.3. There exists a natural bijection between
the set of ${\cal A}$-graded ideals in the polynomial ring $k[{\bf x}]$
and the set of closed points of the scheme ${\cal P}_{\cal A}$.

\noindent {\sl Proof: }
With every point $ f \in {\cal P}_{\cal A} $ we associate the
following binomial ideal in $\, k[{\bf x}] $:
$$ I_f \quad := \quad \langle \,
f^{\bf b}_{\bf u} \cdot {\bf x}^{\bf v}  \,-\,
f^{\bf b}_{\bf v} \cdot {\bf x}^{\bf u} \,\,\,:\,\,\,
{\bf u},{\bf v} \in {\bf N}^n ,\, {\bf b} \in {\bf N}{\cal A}
\,\,\hbox{ s.t. }\,\, deg({\bf u}) = deg({\bf v}) = {\bf b} \,\, \rangle .
\quad \eqno (5.4) $$
Fix ${\bf b} \in  Z_r({\cal A}) \cap {\bf N}{\cal A}$.
There exists ${\bf u} \in deg^{-1}({\bf b})$ such that
$\,f^{\bf b}_{\bf u} \,\not= \, 0$. Every monomial ${\bf x}^{\bf v}$
 of degree ${\bf b}$ is a scalar multiple of ${\bf x}^{\bf u}$
modulo the relations in $I_f$. Therefore
$\,dim_k \bigl( (S/I_f)_{\bf b} \bigr) \leq 1 $.
We must show that equality holds. The equations (5.3) imply
that  the ${\bf b}$-th graded component of $I_f$ is
spanned as a $k$-vector space by the binomials
$\,f^{\bf b}_{\bf u} \cdot {\bf x}^{\bf v}  \,-\,
 f^{\bf b}_{\bf v} \cdot {\bf x}^{\bf u} $
where ${\bf u},{\bf v} \in deg^{-1}({\bf b})$.
This space is a proper subspace of $S_{\bf b}$.
We conclude that $I_f$ is an ${\cal A}$-graded ideal
in $S^{(r)}$, and by Proposition 5.1, it lifts to a unique
${\cal A}$-graded ideal in $S$.

We shall construct the inverse to the map $f \mapsto I_f$.
Let $J$ be any ${\cal A}$-graded ideal in $S$. For each ${\bf b}
\in Z_r({\cal A}) \cap {\bf N}{\cal A}$ there exists a monomial
${\bf x}^{\bf u}$ which does not lie in $J$.
We define $f = f(J)\in {\bf P}$ as follows: for
${\bf v} \in deg^{-1}({\bf b})$
let $f^{\bf b}_{\bf v}$
be the unique scalar in $k$ such ${\bf x}^{\bf v} - f^{\bf b}_{\bf v}
\cdot {\bf x}^{\bf u}$ lies in $J$. Note $f^{\bf b}_{\bf u} = 1$.
We see that $\, f^{\bf b} \, = \, \bigl( \, f^{\bf b}_{\bf v} \,:\,
{\bf v} \in deg^{-1}({\bf b})\bigr)\,$ is a well-defined point
in the projective space ${\bf P}_{\bf b}$, independent of the choice
of ${\bf u}$. Our assumption that $J$ is ${\cal A}$-graded implies
that the two binomials
$$\,{\bf x}^{\bf w} \cdot ( f^{\bf b}_{\bf u}  \cdot {\bf x}^{\bf v}
\, - \, f^{\bf b}_{\bf v}  \cdot {\bf x}^{\bf u} ) \quad \hbox{and} \quad
f^{{\bf b}+{\bf c}}_{{\bf u}+{\bf w}} \cdot  {\bf x}^{{\bf v}+{\bf w}} -
f^{{\bf b}+{\bf c}}_{{\bf v}+{\bf w}} \cdot  {\bf x}^{{\bf u}+{\bf w}} $$
are non-zero constant multiples of each other, whenever
 ${\bf x}^{{\bf u}+{\bf w}}$ is a monomial not in $J$.
This proves that $f$ satisfies the
equations (5.3) and hence lies in ${\cal P}_{\cal A}$.
It is now obvious that $\, I_f \,= \, J$, and we are done. \Box

\proclaim Corollary 5.4. All irreducible components
of ${\cal P}_{\cal A}$ are rational varieties.

\noindent {\sl Proof: }  This follows from the decomposition
theorem for arbitrary binomial schemes in [5]. Indeed,
under the Segre embedding of the product ${\bf P}$, the
equations (5.3) translate into linear equations
with two terms. This shows that ${\cal P}_{\cal A}$ is a binomial scheme.
\Box

\vskip .1cm

There is a natural action  of the torus $(k^*)^n$ on the product
of projective spaces ${\bf P}$. If $f \in {\bf P}$ and
$\lambda \in (k^*)^n$, then $\lambda f$ has coordinates
$\, (\lambda f)^{\bf b}_{\bf u} \,:=\,
\lambda^{\bf u} \cdot f^{\bf b}_{\bf u}$.
Clearly, this action preserves the subscheme ${\cal P}_{\cal A}$,
and we have the relation $ \, I_{\lambda f}\,\, = \,\,
\lambda^{-1} \cdot I_f $.

\proclaim Remark 5.5. The bijection $f \mapsto I_f$ between
 ${\cal P}_{\cal A}$ and ${\cal A}$-graded ideals is $(k^*)^n$-equivariant.

\proclaim Corollary 5.6. The set of isomorphism classes of
${\cal A}$-graded algebras in $k[{\bf x}]$ is in bijection with the set of
$(k^*)^n$-orbits in ${\cal P}_{\cal A}$.

In Section 3 we have seen that for $n=4$ and $d=1$ there may
be infinitely many $(k^*)^n$-orbits on ${\cal P}_{\cal A}$.
In this situation it is desirable to construct a
{\it moduli space} $\, {\cal M}_{\cal A} \, = \,
{\cal P}_{\cal A}/ (k^*)^n \,$ of isomorphism classes
of ${\cal A}$-graded algebras. However, such an enterprise
would immediately face the usual intricacies of
geometric invariant theory [13], such as:
\item{$ \bullet $} The GIT-quotient is not unique but
depends on the  choice of linearization.
Is there a best linearization ? \ (Or perhaps the
Chow quotient of [7] is more useful here~?)
\item{$ \bullet $} Which are the semi-stable orbits~?
And which ones get necessarily lost under the quotient construction~?

\vskip .1cm

\noindent
Leaving these questions for future studies (see Problem 6.5),
we close this section with a corollary about the ``coherent component''
of ${\cal P}_{\cal A}$. Let $e$ denote the point in $\, {\cal P}_{\cal A}
\subset {\bf P} \,$ all of
whose coordinates $e^{\bf b}_{\bf u}$ are equal to one.
Then $I_e $ equals the toric ideal $I_{\cal A}$.

\proclaim Corollary 5.7. The map $f \mapsto I_f$ defines
a bijection between the closure of the $(k^*)^n$-orbit
of $e$ in ${\cal P}_{\cal A}$ and the set of
coherent ${\cal A}$-graded ideals in  $k[{\bf x}]$.
This orbit closure $\,\overline{(k^*)^n \cdot e }$ equals the projective toric
variety defined by the state polytope of ${\cal A}$.

\beginsection 6. Open problems.

Starting from the examples in Section 2, it is easy to construct
incoherent ${\cal A}$-graded algebras for all choices of $n$ and
$d = dim({\cal A})$ with $n \geq d + 3$. The case $n \leq d+1$ being trivial,
and the case $n=3,d=1$ being answered by Theorem 1.1,
the question remains what happens for
$n = d+ 2 \geq 4 $. In view of C.~Lee's result [12] that all
triangulations of  $(d-1)$-polytopes with $d+2$ vertices are coherent,
we venture the following conjecture.

\proclaim Conjecture 6.1.
If $n - d \leq 2 $, then all ${\cal A}$-graded algebras are coherent.

\vskip .1cm

We also conjecture the following converse to Theorem 3.1.

\proclaim Conjecture 6.2. For every polyhedral subdivision
$\Delta$ of a finite set ${\cal A} \subset {\bf N}^d$ there exists an ${\cal
A}$-graded
ideal $I$ whose radical $Rad(I)$ equals the Stanley ideal $I_\Delta$.

Theorem 4.1 and Conjecture 6.2 (if true)  would completely
characterize the reduced schemes defined by ${\cal A}$-graded algebras:
they are precisely the (algebraic sets associated with) polyhedral
subdivisions of ${\cal A}$.
Another obvious question we left open is the following.

\proclaim Problem 6.3. Find an optimal bound for the integer $r$
in Proposition 5.1.

\vskip .1cm

The parameter space ${\cal P}_{\cal A}$ has the following general structure.
It has one nice component consisting of all coherent AGAs
(cf.~Corollary 5.7), and it may have many other components about
which we know very little.
For instance, we do not know whether there
can be embedded components. One line of attack is suggested by the geometry
of Example 3.2. Here the extraneous component
corresponds to a family of incoherent subdivisions of ${\cal A}$.
On the other hand, such incoherent subdivisions
give rise to extraneous components in the
{\it inverse limit of toric GIT-quotients} introduced in [7, \S 4].
Here we make the assumption that all vectors in ${\cal A}$
have the same coordinate sum, so that ${\cal A}$ defines  an
action of the torus $(k^*)^d$ on projective space $P^{n-1} \,$
(see [7] for details).

\proclaim Problem 6.4.
Does there exist a natural morphism from the parameter
space ${\cal P}_{\cal A}$
onto the inverse limit of all toric GIT-quotients
$P^{n-1}/ (k^*)^n$ with respect to the action
$$ (\,x_1 \,:\, x_2 \,: \,\cdots \,: \,x_n\,) \,\,\, \mapsto\,\,\,
( {\bf t}^{{\bf a}_1} x_1 \,:\,{\bf t}^{{\bf a}_2} x_2 \,:\,\,\cdots\,\,
{\bf t}^{{\bf a}_n} x_n \,) \, ? $$

The restriction of the desired morphism to the coherent component is
well-known in combinatorial algebraic geometry: it is the
contraction from the toric variety of the state polytope
onto the toric variety of the secondary polytope [8].
The latter is the Chow quotient $\,P^{n-1}/\!/ \, (k^*)^n\,$ which appears
as the distinguished component in the inverse limit of GIT-quotients.
Our last question was already asked in the end of Section 5.

\proclaim Problem 6.5.
Construct and study the moduli space ${\cal M}_{\cal A}$ of
$(k^*)^n$-orbits on ${\cal P}_{\cal A}$.

\vskip 1.5cm

\noindent
{\sl Acknowledgement: }
Research on this project was supported in part by the
David and Lucile Packard Foundation, the
National Science Foundation, and the Courant Institute
of New York University. I wish to thank David
Eisenbud for many helpful discussions.

\beginsection References

\item{[1]} V.I.~Arnold: ``A-graded algebras and continued fractions'',
{\sl Communications in Pure and Applied Math,} {\bf 42} (1989) 993--1000.

\item{[2]} D.~Bayer and I.~Morrison: ``Gr\"obner bases and geometric
invariant theory I'', {\sl Journal of Symbolic Computation}
{\bf 6} (1988) 209--217.

\item{[3]} B.~Bayer and M.~Stillman: MACAULAY -- A computer algebra
system for algebraic geometry. Available via
anonymous ftp from {\tt zariski.harvard.edu}.

\item{[4]} P.~Diaconis, R.~Graham and B.~Sturmfels:
``Primitive partition identities'', in:
{\sl Paul Erd\"os is 80. Vol II}, Janos Bolyai Society,
Budapest, 1995, pp.~1--20.

\item{[5]} D.~Eisenbud and B.~Sturmfels:
``Binomial ideals'', Manuscript, 1994.

\item{[6]} I.M.~Gel'fand, M.~Kapranov and A.~Zelevinsky:
{\sl Discriminants, Resultants and Multi-Dimensional Determinants},
Birkh\"auser, Boston, 1994.

\item{[7]} M.~Kapranov, B.~Sturmfels and A.~Zelevinsky:
``Quotients of toric varieties'',
{\sl Mathematische Annalen} {\bf 290} (1991), 643-655.

\item{[8]} M.~Kapranov, B.~Sturmfels and A.~Zelevinsky:
``Chow polytopes and general resultants'',
{\it Duke Mathematical Journal} {\bf 67} (1992) 189--218.

\item{[9]} E.~Korkina: ``Classification of A-graded algebras with
3 generators'', {\sl Indagationes Mathematicae} {\bf 3} (1992) 27--40.

\item{[10]} E.~Korkina, G.~Post and M.~Roelofs:
``Alg\'ebre gradu\'ees de type A'', {\sl Comptes Rendues
Acad.~Sci.~Paris, S\'erie I} {\bf 314} (1992) 653--655.

\item{[11]} E.~Korkina, G.~Post and M.~Roelofs:
``Classification of generalized A-graded algebras with 3 generators''
{\sl Bulletin des Sciences Math\'ematiques}, to appear.

\vskip .1cm
\item{[12]} C.~Lee:  ``Regular triangulations of convex polytopes'', in
{\sl Applied Geometry and Discrete Mathematics - The Victor Klee
Festschrift}, (P.~Gritzmann and B.~Sturmfels, Eds.), American
Math.~Soc, DIMACS Series {\bf 4}, Providence, R.I. (1991), 443--456.

\item{[13]} D.~Mumford and J.~Fogarty,
Geometric Invariant Theory (2nd ed.), Springer, 1982.

\item{[14]} H.M.~M\"oller and T.~Mora:
``Upper and lower bounds for the degree of Gr\"obner bases'',
EUROSAM '84, {\sl Springer Lecture Notes in Computer Science}
{\bf 174} (1984) 172--183.

\item{[15]} B.~Sturmfels: ``Gr\"obner bases of toric varieties'',
{\sl T\^ohoku Math.~J.} {\bf 43} (1991) 249-261.

\item{[16]} B.~Sturmfels and R.R.~Thomas: ``Variation of cost functions
in integer programming'', Manuscript, 1994.

\item{[17]} G.~Ziegler: {\sl Lectures on Polytopes},
Graduate Texts in Mathematics, Springer Verlag, New York, 1994.

\bye